\documentclass[11pt]{article}
\usepackage{amsmath,amsthm,amsfonts,graphicx,amssymb}
\usepackage{url}

\usepackage{amsmath}
\usepackage{amsmath}
\usepackage[all]{xy}

\newtheorem{Theorem}{Theorem}[section]
\newtheorem{lemma}[Theorem]{Lemma}
\newtheorem{claim}[Theorem]{Claim}

\theoremstyle{definition}
\newtheorem{definition}[Theorem]{Definition}

\newtheorem{remark}[Theorem]{Remark}

\numberwithin{equation}{section}

\def\ep{\hfill $\Box$}

\def\bp{\noindent{\bf Proof}  \ }

\begin{document}

\begin{center}
{\large {\bf ON THE QUOTIENT-LIFT MATROID RELATION}}
\\[5ex]
 {\bf Jos\'e F. De Jes\'us  (University of Puerto Rico, San Juan,Puerto Rico, U.S.A.) }
\\[4ex]
 {\bf Alexander Kelmans (University of Puerto Rico, San Juan,Puerto Rico, U.S.A.) }

\end{center}

\date{}

\vskip 3ex

\begin{abstract}

It is well known that a matroid $L$ is a lift of a matroid $M$ if and only if every circuit of $L$ is the union of some circuits of $M$. 
In this paper we give a simpler proof of this important theorem.
We also described a discrete homotopy theorem on two matroids of different ranks on the same ground set.


 \vskip 1ex


\end{abstract}

\section{Introduction} 

\label{Intro}

\indent

It is well known that a matroid $L$ is a lift of a matroid $M$ if and only if every circuit of $L$ is the union of some circuits of $M$ \cite{Ox}. 
The proof of this characterization given in the classic book "Matroid Theory", by James Oaxley, depends heavily on mathematical induction, assumes the elementary quotient construction and is based on the notion of modular cuts of flats.
In this paper we give a simpler (constructive) proof of this important theorem.
Our proof is based only on the simple notions of cyclic sets and matroid ciclomaticity.
We also prove a discrete homotopy theorem on two matroids $M_1$ and $M_2$ of different ranks $r_1 < r_2$ on the same ground set $E$ saying that $M_2$ can be obtained from $M_1$ by a series of $r_2 - r_1$ one element extentions.

 \section{Preliminaries}

All notions and basic facts on matroids that are used here can be found in   \cite{Ox,Wlsh}.
In this section we will remind the reader the main matroid notions and facts we need for our proof.
 
Given a family ${\cal F}$ of subsets of a set $E$ (i.e. ${\cal F}\subseteq 2^E$), let 
${\cal M}in \hskip 0.5ex {\cal F}$ 
and  ${\cal M}ax \hskip 0.5ex {\cal F}$ denote
 the family  of elements of $S$ which are minimal and maximal by the set inclusion,
 respectively. 
 A family ${\cal P}$ of subsets of  $E$ is called a  {\em  clutter} if 
 $P, Q \in  {\cal P}~and ~ P \subseteq Q \Rightarrow P = Q$, 
 and so
 ${\cal M}in \hskip 0.5ex {\cal F}$ and ${\cal M}ax \hskip 0.5ex {\cal F}$ are clutters.

\vskip 1.0ex

Let $E$ be a finite non-empty set and $\cal{I}$ a family of subsets of $E$, i.e. 
${\cal I} \subseteq 2^E$.
A pair $M = (E, {\cal I})$ is called a {\em  matroid} if
\vskip 0.7ex
\noindent
$(AI0)$ $\emptyset \in  {\cal I}$,
\vskip 0.7ex
\noindent
$(AI1)$ if $X \in  {\cal I}$ and $Z  \subseteq X$, then $Z \in  {\cal I}$, and
\vskip 0.7ex
\noindent
$(AI2)$ if $X,Y \in  {\cal I}$ and  $|X| < |Y|$, then there exists   
$y  \in Y\setminus X$ such that $X + y \in  {\cal I}$.

\vskip 1.0ex

The set  $E$ is called the  {\em  ground set of} $M$ and an element of  ${\cal I}$ is called an {\em independent set} of $M$.  
The family ${\cal I}  = {\cal I}(M)$ is   the   family of {\em independent sets} of $M$,
${\cal D} =  {\cal D}(M) = 2^E \setminus \cal{I}$ is the family of
 {\em dependent sets} of $M$,
${\cal B} = {\cal B}(M)= {\cal M}ax \hskip 0.5ex \cal{I}$ is the family of {\em  bases} of $M$,  and
${\cal C} = {\cal C}(M) ={\cal M}in \hskip 0.5ex \cal{D}$ is the family of {\em circuits} of $M$. 

Let ${\cal B}^* = {\cal B^*}(M) = \{E \setminus B: B \in {\cal B}(M) \}$. It is easy to see that ${\cal B}^*$ is the set of bases of a matroid (denoted by $M^*$) on the ground set $E$. Matroids $M$ and $M^*$ are called  {\em dual matroids}.

\vskip 1ex

Given a matroid $M = (E, {\cal I})$ and $Z \subset E$, let $E' = E \setminus Z$ and 
$I' = \{I  \in  {\cal I}:   I \subseteq E' \}$. Then, obviously, $R= (E', I')$ is a matroid.
Put $R = M \setminus Z$. We say that 

$(d)$ {\em  $R = M \setminus Z$ is obtained from $M$ by deleting set $Z$} and  

$(c)$ {\em  $R^*$ is obtained from $M$ by codeleting $($or contracting$)$ 
set $Z$ and put $R^* =  M / Z$}.

\vskip 1ex

Let $\rho(M) =  |B|$, where $B \in  {\cal B} (M)$.

\vskip 1.5ex

Obviously, the family ${\cal C}= {\cal C}(M)$ of circuits of a matroid $M$ has the following properties:
\vskip 0.7ex
\noindent
$(AC0)$ $\emptyset \not \in  {\cal C}$ and
\vskip 0.7ex
\noindent
$(AC1)$ $C_1, C_2 \in  {\cal C}~and ~ C_1\subseteq C_2 \Rightarrow C_1 = C_2$,
and so ${\cal C}$ is a clutter.

Moreover, ${\cal C}$ is the set of circuits of a matroid if and only if
${\cal C}$ satisfies axioms $(AC1)$, $(AC2)$, and the following axiom
\vskip 0.7ex
\noindent
$(AC2)$ if $C_1, C_2 \in  {\cal C}$,  $C_1\ne C_2$, and $e \in C_1\cap C_2$, then
there exists $C \in  {\cal C}$ such that $C \subseteq C_1\cup C_2 - e$.
\vskip 0.7ex
\noindent
Axiom $(AC2)$ is called the  {\em circuit elimination axiom} of a matroid
(CEA, for short).

\vskip 1ex

It turns out   that the set ${\cal C}$  of circuits of a matroid also satisfies the following
{\em strong circuit elimination axiom} of a matroid (SCEA, for short):
\vskip 0.7ex
\noindent
$(AC2!)$ if $C_1, C_2 \in  {\cal C}$,  $C_1\ne C_2$,  
$e \in C_1\cap C_2$ and $d \in C_1\setminus C_2$, then there
exists $C \in  {\cal C}$ such that $d \in C \subseteq C_1\cup C_2 - e$.

\vskip 1ex

Consider  an independent
 set $I$ of a matroid 
$M = (E, {\cal I})$, $x \in E \setminus I$. If  $I + x$ is not independent, then
by $(AC2)$
there  exists a unique circuit $C$ of $M$ such that $C \subseteq I + x$ and, obviously, $x \in C$. We call $C$ the {\em  $x$-fundamental circuit of $M$}
({\em  with respect to  $I$}) and denote it $C(x,I)$.

 \vskip 1ex

We call a subset $A$ of $E$
a {\em cyclic set} of $M$ if $A$ is the union of some circuits of $M$. 

\vskip 1ex

Let $E$ and $X$ be non-empty finite disjoint sets. 
Let  $N$ be a matroid on the ground set $E \cup X$ and ${\cal C}_N$ 
the set of circuits of $N$.

As above, $N \setminus X$ is the matroid on $E$ obtained from $N$ by {\em deleting} set $X$ and $N / X$ is  the matroid on $E$ obtained from N by {\em contracting} set $X$. Obviously, the circuits of 
$N \setminus X$ are the circuits of ${\cal C}_N$ that are contained in $E$.
Given two sets $Y$ and $Z$, we call $Y \cap Z$ the {\em trace of $Z$ in $Y$} and also  the {\em trace of $Y$ in $Z$}.
Obviously,  $R$ is a circuit of $M = N / X$ if and only if $R$ is a minimal trace of a  circuit of $N$ in $E$.

\section{Main results} 

Let $M$ be   a matroid on a finite set $E$.

\begin{definition} 
\label{cdpair}

 Let $M$ and $L$ be matroids on $E$. 
 We call $(M, L)$ an
 {\em $X$-codeletion-deletion pair} (or simply, 
 an {\em  $X$-$(c,d) $-pair} or just a {\em $(c,d) $-pair})
 if the exists  a finite non-empty set $X$ disjoint from $E$ and  a matroid $N$ on 
 $E  \cup X$ such that  $M = N / X$ and 
 $L = M \setminus X$. In this case $M$ is also called a  {\em  quotient of} $L$ and
 $L$ is called a  {\em  lift of} $M$, and so we can also call $(M, L)$ 
 a {\em quotient-lift  pair}.

  \end{definition}

 \begin{Theorem} 
\label{liftproperty}

Let $M$ and $L$ be matroids on $E$.  Suppose that $(M, L)$ is an 
$X$-$(c,d) $-pair for some set $X$. 
Then every circuit of $L$ is the union of some circuits of $M$.

\end{Theorem}
 
\bp Since $(M, L)$ is an $X$-$(c,d) $-pair, we have:  $X \cap E = \emptyset $ and there exists matroid $N$ on $X \cup E$ such that 
$N / X = M$ and   $N \setminus X = L$. 
Also  ${\cal C}(L)\subseteq {\cal C}(N)$.

Let $D \in {\cal C}(L)$.  If $D \in {\cal C}(M)$, then we are done. 
So suppose that  $D \not \in {\cal C}(M)$. 
Since $(M, L)$ is an $X$-$(c,d) $-pair,
there exists $C \in {\cal C}(M)$ such that $C \subset D$ and $C$ is a minimal trace of a circuit of $N$, say $Q$, in $E$. We need to prove that every element $d$ of $D$ is in some minimal trace of a circuit of $N$ in $E$ which is a subset of $D$. If $d \in C$, we are done. So suppose $d \not \in C$.

Note that $Q \cap X \ne \emptyset$  for otherwise 
$C \in {\cal C}(L)$, a contradiction. Since $C \ne \emptyset$ and $C \subset Q \cap D$, there exits $a \in  Q \cap D$. 
By (SCEA),  applied to circuits $Q$ and $D$ in $N$ with 
$a \in Q \cap D$ and $d \in D \setminus Q$, there exists a circuit $S$ of $N$ such that 
$d  \in S \subseteq (Q \cup D) - a$, and so $d$ is 
in the trace $S'$ of $S$ in $E$.
Let ${\cal T}'$ be the set of all traces of circuits of $N$ in $E$ containing $d$ which are subsets of $D$. Let  $T'$ be a minimal set in ${\cal T}'$.
If $T'$ is a minimal trace of a circuit of $N$  in $E$, then we are done. So suppose that $T'$ is not a minimal trace  of a circuit  of $N$ in $E$. 
Then there is  a circuit $Z$ of $N$ with  trace $Z'$ in $E$ and $z \in Z'$ such that $Z'  \subseteq T' - d \subseteq D$. 
Now 
by (SCEA),  applied to circuits $Z$ and $T$ in $N$ with 
$z \in Z \cap T$ and $d \in T \setminus Z$, there exists a circuit $P$ of $N$ such that 
$d  \in P \subseteq (Z \cup T) - z$. Let  $P'$ be the trace of $P$ in $E$.
Then $d  \in P' \subseteq (Z' \cup T') - z  \subseteq  T' - z$.
Thus, $T'$ is not a minimal trace  of a circuit in $N$ containing  $d$, 
a contradiction.
 \ep

\vskip 1.5ex

Below we will show that the converse of Theorem \ref{liftproperty} is also true. We need some more definitions and preliminary facts.

\begin{definition} $(${\sc Fundamental $s$-family of circuits in a matroid}$)$ 
 \label{FundList}
\\
Let ${\cal F} \subseteq {\cal C}(M)$. We call  ${\cal F}$ a {\em fundamental 
$s$-family of circuits} of $M$ {\em with respect to $I$} if there exists  $I \in {\cal I}(M)$ and 
$S \subseteq  E \setminus I$ such that ${\cal F} = \{C(x,I): x \in S\}$ and 
$|S| = s$. 
\end{definition}

\begin{definition} 
\label{cyclomaticity}
Let $A \subseteq E$ and  $I$  be  a maximal independent subset of $A$.
We call $c(A) = | A \setminus I|$   {\em the cyclomaticity of} $A$ (also known as  {\em the cyclomatic number} or {\em the nullity of} $A$).

\end{definition}

We will need the following well-known fact due to J. Edmonds.

\begin{lemma} $(${\sc Spanning property of a fundamental family of circuits  in a matroid}$)$
\label{spanprop}
Suppose that $A \subseteq E$, $A$ is a dependent set of $M$,   $I$ is a maximal independent set of $A$ in $M$, and $c = c(A) = |A \setminus I|$. Then 

$~~~~~~~~~~~\cup \{C \in {\cal C}(M): C \subseteq A\} = \cup \{C(x,I): x \in A \setminus I\}$ 
\\
and 
$\{C(x,I): x \in A \setminus I\} $
is a fundamental $c$-family  of circuits of $M$ with respect to $I$.
\end{lemma}

\begin{claim} 
\label{A,Q}

Suppose that $A$ is a cyclic set of $M$, 
$I$ is a maximal independent subset of $A$,  
${\cal F} = \{C(a, I): a \in A \setminus I\}$ is a fundamental $c$-family  of circuits 
$ C(a, I)$ of $M$ such that $\cup \{F \in {\cal F}\} = A$, and (as above) 
$|A \setminus I| = c(A) = c$.
Then for every circuit $Q$ of $M$ which is not a subset of $A$  there exists a circuit $Q'$ of $M$ such that 
$Q'$ is a subset of $A \cup Q$
and
${\cal F}' = {\cal F} \cup \{Q'\}$ is a fundamental $(c +1)$-list 
of circuits of $M$.
\end{claim}

\bp  Let $A' = A \cup Q$ and $I'$ be a maximal independent set in $  A'$
such that $I  \subseteq I'$. Then every $C(a, I)$ is also a fundamental circuit (rooted at $a$) with respect to the independent set $I'$.
Obviously, $Q  \setminus I' \ne \emptyset $, say $q \in Q  \setminus I'$, and 
$q \not \in A \setminus I$.
Then $I' + q$ has a unique circuit $Q' = C(q,I')$ of $M$ containing $q$ and 
$Q' \subset A \cup Q$. Thus, ${\cal F}' = {\cal F} \cup \{Q'\}$ is a fundamental 
$(c +1)$-list  of circuits of $M$. 
\ep

\vskip 1.5ex

From Claim \ref{A,Q} we have:  

\begin{lemma} $(${\sc Extension  property of cyclic sets in a matroid}$)$
\label{cycl-sets-extnsion}
\\
Let $A_1$ and $A_2$ be distinct cyclic sets of $M$ such that $A_1 \not \subseteq A_2$ and let $c(A_1) = c$. 
Then there exists a cyclic set $A$ of $M$ such that $A \subseteq A_1 \cup A_2$ and $c(A) = c+1$.
\end{lemma}

Using Lemma \ref{spanprop}, it is easy to prove the following  

\begin{claim} 
\label{A,B,D}

Let $A$ be a cyclic set of $M$, $c(A) = c$, $D \in {\cal C}(M)$, 
and $d \in D \subseteq A$. Then there exists a base $B$ of $M$ and a fundamental $c$-family ${\cal F}$  of circuits of $M$ $($with respect to $B$$)$ such  that $D  \in {\cal F}$, $\cup\{F \in {\cal F}\} = A$, and $d \notin B$.

\end{claim}

\bp Obviously, $D - d$ is an independent set of $M$. Let $I$ be a maximal independent set in $A$ such that   $D - d  \subseteq I$. 
By Lemma   \ref{spanprop},

$~~~~~~~~\cup \{C \in {\cal C}(M): C \subseteq A\} = \cup \{C(x,I): x \in A \setminus I\}$.
\\
Since $D - d  \subseteq I$, clearly $C(d,I) = D$. Let $B$ be a base of $M$ containing $I$ as a subset and ${\cal F} = \{C(x,I): x \in A \setminus I\}$. Then 
${\cal F}$ is a fundamental $c$-family 
 of circuits of $M$ $($with respect to $B$$)$, 
 $D\in {\cal F}$, 
 $\cup \{F  \in {\cal F}\} = A$, and $d \notin B$.
\ep

\begin{claim} 
\label{cycl-sets-elimination1}

Let $A$  be a cyclic set of $M$ with $c(A) = c$ and $Q$ be a circuit of $M$ which is not a subset of $A$. 
Then for every $a \in A \cap Q$ there exists a cyclic set $A'$ of $M$ such that $a \notin A'  \subseteq A \cup Q$ and $c(A') = c$.

\end{claim}

\bp By Lemma \ref{spanprop},
$A = \cup \{C \in {\cal C}(M): C \subseteq A\} = \cup \{C(x,I): x \in A \setminus I\}$, 
where  $c =  |A \setminus I|$. 
Let $I'$ be a maximal independent set in $  A \cup Q$ such that $I  \subseteq I'$.
Then for every $q \in Q \setminus I'$ there is a unique circuit $Z$ such that
$Z$ is a fundamental circuit $C(q, I')$ with respect $I'$ rooted at $q$.

First, suppose that $a \in A \setminus I$. 
Let  $A' = (A  \setminus C(a, A))  \cup  C(q,I')$. 
Then $A'$ is a required set.

Now suppose that $a \in I$. Then 
$a \in C(z,I) - z$ for some $z \in A\setminus I$. Also by (SCEA) in $M$, there exists a circuit $Q'$ of $M$ such that $a \notin Q'$ and $q \in Q'
  \subseteq Q \cup C(z,I)$. Let  
  $A' = (A  \setminus C(z, A))  \cup  C(q,I')$.
  Then again $A'$ is a required set.
\ep

\vskip 1.5ex

From Claim \ref{cycl-sets-elimination1} we have:  

\begin{lemma} $(${\sc Elimination  property of cyclic sets in a matroid}$)$
\label{cycl-sets-elimination2}
\\
Let $A_1$ and $A_2$ be distinct cyclic sets of $M$ and let $c(A_1) = c$. 
Then for every $a \in A_1 \cap A_2$ there exists a cyclic set $A$ of $M$ such that $a \notin A  \subseteq A_1 \cup A_2$ and $c(A) = c$.
\end{lemma}

\begin{claim} 
\label{r(L)>r(M)}

Let $M$ and $L$ be distinct matroids on $E$.  Suppose that $(M, L)$ is an 
$X$-$(c,d) $-pair for some set $X$. Then $\rho(L) > \rho(M)$. 

\end{claim}

\bp By Theorem \ref{liftproperty}, every circuit of $L$ is the union of some circuits of $M$. Therefore every dependent set of $L$ is also a dependent set of $M$ or, equivalently, every independent set of $M$ is  an independent set of $L$. In particular, every base of $M$ is an independent set of $L$. Therefore $\rho(L) \ge \rho(M)$. Since $M \ne L$, clearly 
$B \in {\cal B}(M) \Rightarrow B \in {\cal I}(L)  \setminus {\cal B}(L)$.
Thus,  $\rho(L) > \rho(M)$.
\ep

\vskip 1.5ex 

We need the following  
 
\begin{lemma} 
\label{bigcyclo}

Let $M$ and $L$ be matroids on $E$ and $\rho(L) = \rho(M) + s$, where 
$s \in \mathbb{N}$. Suppose that every circuit of $L$ is the union of some circuits of $M$. If $A$ is a cyclic set of $M$ with $c(A) = s+1$, then $A$ is not an independent set of $L$.
\end{lemma} 

\bp Suppose, on the contrary, that  $A$ is an independent set of $L$. Since $A$ is a cyclic set of $M$ with  $c(A) = s+1$,
there exists a subset $R$ of $A$ with $s+1$  elements such that $I = A  \setminus R$ is a maximal independent set of $A$ in $M$.   Let $B$ be a base of $M$ 
such that $I \subseteq B$. 

Note that $|R \cup B| = \rho(M) + s + 1 > \rho(L)$.  Thus, $R \cup B$ is a dependent set of $L$, i.e. there exists  $D  \in {\cal C}(L)$ such that $D \subseteq R \cup B$. Since $A$ is an independent set of $L$, clearly 
$D \not \subseteq A$. 
Let $d \in D \setminus A$. Since $D$ is the union of some circuits of $M$, there exists  $D'  \in {\cal C}(M)$ such that $d \in D' \subseteq D \subseteq R \cup B$.  Now $\cup \{C(x,I): x \in A \setminus I\} = \cup \{C(x,I): x \in R\} = \cup \{C(x,B): x \in R\}$. By Lemma \ref{spanprop}, $ \cup \{C(x,B): x \in R\}  = \cup \{C \in {\cal C}(M): C \subseteq B \cup R \}$. It follows that $D' \in \cup \{C(x,I): x \in A \setminus I\}$, and so $d \in D' \subseteq A$. However $d \in  D \setminus A$, a contradiction. 
 \ep

\vskip 1.5ex 

Now we are ready to prove the converse of Theorem \ref{liftproperty}. 
By Claim \ref{r(L)>r(M)}, 
if $(M, L)$ is an  $X$-$(c,d) $-pair of matroids for some set $X$, then 
$\rho(L) - \rho(M) = s \in \mathbb{N}$.
Therefore in the converse  of Theorem \ref{liftproperty} we can assume that 
$\rho(L) - \rho(M) = s \in \mathbb{N}$.

\begin{Theorem} 
\label{liftcriterion}
Let $M$ and $L$ be matroids on $E$ and $\rho(L) = \rho(M) + s$, where 
$s \in \mathbb{N}$. Suppose that every circuit of $L$ is the union of some circuits of $M$. Then $(M, L)$ is an 
$X$-$(c,d) $-pair for some set $X$ with $s$ elements. 

\end{Theorem}

\bp Let $X$ be a set with $s$ elements.
By definition, a pair $(M, L)$ is an  $X$-$(c,d) $-pair if and only if  
there exists a matroid $N$ on $E \cup X$ such that  $N / X = M$ and  $N \setminus X = L$.

Let 

\vskip 0.5ex

$~~~~~~~~~{\cal X} = \cup \{A \cup Z : A \in {\cal CS}(M) \cap {\cal I}(L), \emptyset \ne Z  \subseteq X, 
 ~and~c(A) + |Z| = s + 1\}$.

We prove that 
${\cal C}(L) \cup {\cal X}$ satisfies the elimination axiom of the set of circuits of a matroid, say  $N$, on $E \cup X$, i.e.  that ${\cal C}(L) \cup {\cal X}  =  {\cal C}(N)$.
 Let $C_1, C_2 \in {\cal C}(N)$, $C_1 \ne C_2$, and $a \in C_1 \cap C_2$. 

\vskip 1.5ex
\noindent
${\bf (p1)}$ 
Suppose that $C_1,C_2 \in {\cal C}(L)$. Then, obviously, our claim is true.

\vskip 1.5ex
\noindent
${\bf (p2)}$
Suppose that $C_1\in {\cal C}(L)$ and $C_2 \in {\cal X}$. 
Then  
$C_2 = A \cup Z$, where 

$A \in {\cal CS}(M) \cap {\cal I}(L), Z  \subset X, ~and~c(A) + |Z| = s + 1$.
\\
 Since $a \in C_1 \subseteq E $, $Z \subseteq X$, and $X \cap E = \emptyset$, clearly $a \notin Z$, and so $a \in A$. 
Both $A$ and $C_1$   are cyclic sets of $M$. By Lemma \ref{cycl-sets-elimination2} with $A_1=A$ and $A_2=C_1$, there exists a cyclic set $A'$ of $M$ such that $a \notin A'  \subseteq C_1 \cup A \subseteq C_1 \cup C_2$ and $c(A') = c(A)$. If $A'$ contains no circuit of $L$, then $A' \cup Z \in {\cal X}$. Since  $a \notin A'  \cup Z \subseteq  C_1 \cup C_2$, we are done. If $A'$ contains a circuit $C'$ of $L$, then $a \notin C' \subseteq A'  \cup Z \subseteq  C_1 \cup C_2$, and we are also done.

\vskip 1.5ex
\noindent  
${\bf (p3)}$
Suppose that  $C_1,C_2 \in {\cal X}$, namely, $C_1= A_1 \cup Z_1 \in {\cal X}$ and $C_2= A_2 \cup Z_2 \in {\cal X}$.

\vskip 1.5ex
 
${\bf (p3.1)}$
Suppose  that  
$Z_1 = Z_2 = \{a\}$. 
Then $A_1 \ne A_2$ and $c(A_1) = c(A_2) = s$. 
Since $A_1$ and $A_2$ are distinct cyclic sets of $M$, by Lemma \ref{cycl-sets-extnsion},  there exists a cyclic set $A$ of $M$ such that $A \subseteq A_1 \cup A_2$ and $c(A) = c(A_1)+1 = s + 1$. 
Now by Lemma \ref{bigcyclo},
$A$ contains a circuit, say $Q$, of $L$ as a subset, and we are done 
because $Q \subseteq A_1  \cup A_2$ and $a \not \in A_1  \cup A_2$.

\vskip 1.5ex
 
${\bf (p3.2)}$
Now suppose that at least one of  $Z_i$'s, say $Z_1$, has an element distinct from $a$. 
We remind that $(A_1 \cup A_2) \cap (Z_1 \cup Z_2) = \emptyset$. Thus, 
either $a \in Z_1 \cap Z_2$ or $a \in A_1 \cap A_2$. 

  First, suppose that $a \in Z_1 \cap Z_2$. Since $A_1$ and $A_2$ are distinct cyclic sets of $M$, by Lemma \ref{cycl-sets-extnsion},  there exists a cyclic set $A$ of $M$ such that $A \subseteq A_1 \cup A_2$ and $c(A) = c(A_1)+1$. If $A$ contains a circuit of $L$, we are done. If $c(A) = s + 1$, then we are done by the  Lemma \ref{bigcyclo}.  If $c(A) < s +1$, then $c(A_1) < s$, and therefore $|Z_1| > 1$. If $A$ contains no circuit of $L$, then $(A \cup Z_1 \setminus a ) \in {\cal X}$, and we are also done.  

Finally, suppose that $a \in A_1 \cap A_2$. Then by  
Lemma \ref{cycl-sets-elimination2},  there exists a cyclic set $A$ of $M$ such that $a \notin A  \subseteq A_1 \cup A_2$ and $c(A) = c$. 
If $A$ contains a circuit of $L$ as a subset, then we are done. 
If $A$ contains  no circuit of $L$ as a subset, then $A  \cup Z_1 \in {\cal X}$ and again we are done.
 \ep

\vskip 1.5ex

 Using Theorems \ref{liftproperty} and \ref{liftcriterion} it is easy to prove the  following  useful fact.

 \begin{lemma}  {\sc (Transitivity of (c,d)-pair relation between matroids)}
\label{(c,d)-pair-transitivity}
\\
Let  $K$, $L$, and $M$ be matroids on $E$. If $(M,L)$ is an $X$-$($c,d$)$-pair, $(L,K)$ is a $Y$-$($c,d$)$-pair, 
 and $X \cap Y = \emptyset $, then $(M,K)$ is   an $X \cup Y$-$($c,d$)$-pair.
\end{lemma}

\bp By Theorem \ref{liftproperty}, every circuit of $L$ is the union of some circuits of $M$ and every circuit of $K$ is the union of some circuits of $L$.
Therefore every circuit of $K$ is  the union of some circuits of $K$. Thus, by 
 Theorem \ref{liftcriterion},  $(M, K)$ is a $Z$-$(c,d)$-pair for some set $Z$.
 Since  $(M,L)$ is an $X$-$(c,d)$-pair and  $(L,K)$ is a $Y$-$(c,d)$-pair
 (where $X \cap Y = \emptyset $), we have: $Z = X \cup Y$, and so 
 $(M,K)$ is   an $X \cup Y$-$(c,d)$-pair.
\ep

\vskip 1.5ex

Here is another criterion for an $X$-{\em $(c,d)$}-pair of matroids.
                        
 \begin{Theorem} 
\label{delete-contract-pair}

Let $M$ and $L$ be matroids on $E$, $X = \{x_1, \ldots , x_k\}$, and 
$E \cap X = \emptyset $.
Then the following are equivalent:
\vskip 1ex
\noindent
$(c1)$ $(M, L)$ is an $X$-$(c,d) $-pair and
\vskip 1ex
\noindent
$(c2)$ there exists a sequence $(L_0, L_1, \ldots  , L_k)$  such that 
$L_0 = M$, $L_k = L$, 
each $L_i$ is a matroid on $E$, and 
each $(L_{i-1}, L_i)$ is an $x_i$-$(c,d) $-pair, where $1\le i \le k$.
\end{Theorem}
 
 \bp   
Claim  $(c2) \Rightarrow (c1)$ can be easily proved by induction on $|X| = k$ using 
Lemma \ref{(c,d)-pair-transitivity}.

Now we prove Claim $(c1) \Rightarrow (c2)$  by induction on $|X| = k$. 
If $|X| = 1$, then  Claim $(c1) \Rightarrow (c2)$ is obviously true.  Suppose that Claim $(c1) \Rightarrow (c2)$ is true for $|X| = k-1$. We need to prove that  Claim $(c1) \Rightarrow (c2)$ is also true for $|X| = k$. 
Let $X' = X  -  x_k$, and so $|X'| = k-1$. By $(c1)$, there exists a matroid $N$ on  $E \cup X$ with $X = \{x_1, \ldots , x_k\}$  such that 
$M = N/ X = (N / x_k) / X'$ and 
$L =  N \setminus X = (N \setminus x_k) \setminus X'$. 
Let $N' = N / x_k$. 
 Then $N'$ is a matroid on $E \cup X'$ with $X' =\{x_1, \ldots , x_{k-1}\}$.
Let $L' = N'  \setminus X' = (N / x_k) \setminus X'$.
Obviously, $(M, L')$ is an $X'$-$(c,d) $-pair. Put $M = L_0$ and $L' = L_{k-1}$. 
By the induction hypothesis, $(c2)$ holds for pair $(M, L')$, namely,
there exists a sequence $(L_0, L_1, \ldots  , L_{k-1})$  such that 
$L_0 = M$, $L_{k-1} = L'$, 
each $L_i$ is a matroid on $E$, and 
each $(L_{i-1}, L_i)$ is an $x_i$-$(c,d) $-pair, where $1\le i \le k-1$.
Then 
\vskip 1ex
$
~~~~~~~~~~~~~
(N \setminus X')/ x_k = (N / x_k) \setminus X' = L' = L_{k-1}$. 
\vskip 1ex
\noindent
Put $L_k = L$. Then 
$L_k =  (N \setminus x_k) \setminus X' = (N \setminus X')\setminus x_k$. 
Thus, $(L_{k-1}, L_k)$ is an 
$x_k$-$(c,d)$-pair, and so $(c2)$ also holds for $|X| = k$.
 \ep

 \begin{remark} $(${\sc Construction of intermediate matroids in homotopy}$)$
\label{c1-c2}

Claim $(c1) \Rightarrow (c2)$ in Theorem \ref{delete-contract-pair} can also be proved by putting 
\vskip 1ex
${\cal C}(L_i) = \{C \in {\cal C}(L): c(C) \leq i\}    \cup \{A \in {\cal CS}(M):   A \in {\cal I}(L) ~and~ c(A) = i\}$ 
\vskip 1ex
\noindent
and using the arguments  similar to those in the proof of Theorem \ref{liftcriterion}.  

 \end{remark}  

\vskip 1ex

 \begin{remark} 
\label{reduction}

From Theorem \ref{delete-contract-pair}  it follows  that the problem of constructing for a given matroid $M$ all matroids $L$ such that $(M, L)$ is an $X$-$(c,d) $-pair can be reduced to the same problem for $|X| = 1$, i.e. to the problem of constructing all so-called  {\em elementary $(c,d)$-pairs} 
$(M, L)$. 
\end{remark}

\end{document}